\documentstyle[12pt]{article}

\def\ca{{\cal A}}
\def\Z{{\bf Z}}
\def\C{{\bf C}}
\def\cm{{\cal M}}
\def\De{\Delta}
\def\om{\omega}
\def\al{\alpha}
\def\ba{\begin{array}}
\def\ea{\end{array}}

\title{{\Large\bf CHARACTERIZATIONS OF SOME CLASS OF FORMAL POWER SERIES
VIA M\"OBIUS CATEGORIES OF FULL BINOMIAL TYPE}}
\author{EMIL DANIEL SCHWAB}
\date{}

\begin{document}
\maketitle

\begin{quote}
{\small{\bf Abstract.} In [7], [8] and [9] we have obtained characterizations
(we called them characterizations of Lambek-Carlitz type)
of some class of formal power series. Now, we discuss, investigate and
compare these results from a categorial point of view.}
\end{quote}

\begin{center}
{\bf 1.}
\end{center}

J. Lambek [4] proved that an arithmetical function $f$ is completely multiplicative if and
only if it distributes over every Dirichlet convolution $(f(g*_Dh)=fg*_Dfh$ for all arithmetical functions $g$ and $h$). The
Dirichlet convolution of two arithmetical functions $f$ and $g$ is defined by:
$$(f*_Dg)(n)=\sum_{d|n}f(d)g\left(\frac{n}{d}\right)\quad (n,d\in\Z^+).$$

The problem of Carlitz [1] shows that the arithmetical function $f$ is ne\-ce\-ssarily completely
multiplicative if and only if it distributes over the parti\-cu\-lar Dirichlet convolution $\zeta*_D\zeta$, where $\zeta$ is
defined by: $\zeta(n)=1$, $\forall n\in\Z^+$. Thus we have:

{\bf Theorem 1.1.} ([4], [1], the Lambek-Carlitz characterization of completely multiplicative arithmetical functions.) {\it The following
statements are equivalent:

{\rm(1)} $f$ is completely multiplicative;

{\rm (2)} $f(g*_Dh)=fg*_Dfh$ for all arithmetical functions $g$ and $h$;

{\rm (3)} $f(g*_Dg)=fg*_Dfg$ for all arithmetical function $g$;

{\rm(4)} $f\tau=f*_D f$,

\noindent where $\tau(n)$ is the number of positive divisors of
$n$ $(\tau=\zeta*_D\zeta)$.}

The additive version of Theorem 1.1. is the following.

{\bf Theorem 1.2.} {\it The following statements are equivalent:

{\rm (1)} $f$ is completely additive;

{\rm (2)} $f(g*_Dh)=fg*_Dh+g*_Dfh$ for all arithmetical functions $g$ and $h$;

{\rm (3)} $f(g*_Dg)=2(fg*_Dg)$ for all arithmetical function $g$;

{\rm(4)} $f\tau=2(f*_D\zeta)$.}

The unitary convolution of two arithmetical functions $f$ and $g$ is defined by:
$$(f*_U g)(n)=\sum_{d|n,\left(d,\frac{n}{d}\right)=1} f(d)g\left(\frac{n}{d}\right)\quad (n,d\in\Z^+).$$
Let $\ca$ denote the set of all arithmetical functions. Then $(\ca,+,*_U)$ is a ring called the unitary
ring of arithmetical functions. The map
$\eta:\C[[X]]\to\ca$ defined by:
$$\eta\left(\sum_{n=0}^\infty a_nX^n\right)(m)=\om(m)!a_{\om(m)}\quad (m\in\Z^+),$$
where $\om(m)$ denotes the number of distinct prime factors of $m$, is an embedding of the ring of formal power series $\C[[X]]$ in the
unitary ring of arithmetical functions $\ca$
(see [7]). In [7] we have obtained the following characterization of the Lambek-Carlitz type of exponential series:

{\bf Theorem 1.3.} ([7]). {\it Let $\sum_{n=0}^\infty a_n X^n\in\C[[X]]$ be a formal power series
with $a_1\not=0$. The following  statements are equivalent:

{\rm(1)} $\sum_{n=0}^\infty a_n X^n$ is an exponential series $(a_n=\frac{a_1^n}{n!}$, $\forall
n\in\Z^+\cup\{0\})$;

{\rm(2)} $\eta\left(\sum_{n=0}^\infty a_n X^n\right)$ is a multiplicative arithmetical function;

{\rm (3)} $\sum_{n=0}^\infty a_nX^n\odot\left(\sum_{n=0}^\infty b_n X^n\cdot\sum_{n=0}^\infty
c_nX^n\right)=$

\hspace{0.6cm} $=\left(\sum_{n=0}^\infty a_n X^n\odot\sum_{n=0}^\infty b_nX^n\right)\cdot
\left(\sum_{n=0}^\infty a_n X^n\odot
\sum_{n=0}^\infty c_nX^n\right)$,

\hspace{0.6cm} $\forall \sum_{n=0}^\infty b_n X^n$,
$\sum_{n=0}^\infty c_n X^n\in\C[[X]]$;

{\rm(4)} $\sum_{n=0}^\infty a_nX^n\odot\left(\sum_{n=0}^\infty b_n X^n\cdot
\sum_{n=0}^\infty b_n X^n\right)=$

\hspace{0.6cm} $=\left(\sum_{n=0}^\infty a_nX^n\odot\sum_{n=0}^\infty b_n X^n\right)\cdot
\left(\sum_{n=0}^\infty a_n X^n\odot \sum_{n=0}^\infty b_n X^n\right)$,

\hspace{0.6cm} $\forall\sum_{n=0}^\infty b_nX^n\in\C[[X]]$;

{\rm(5)} $\sum_{n=0}^\infty 2^n a_n X^n=\sum_{n=0}^\infty a_n X^n\cdot\sum_{n=0}^\infty a_n X^n$,

\noindent where the ``multiplication" $\odot$ of two formal power
series $\sum_{n=0}^\infty a_n X^n$ and $\sum_{n=0}^\infty b_n X^n$ is defined by:}
$$\sum_{n=0}^\infty a_n X^n\odot \sum_{n=0}^\infty b_n X^n=\sum_{n=0}^\infty n!a_nb_nX^n.$$

(3) is a characterization of Lambek type (distributivity over the product of series) and (5) is a
distributivity condition over a particular
product of series:
$\sum_{n=0}^\infty
\frac{1}{n!}X^n\cdot\sum_{n=0}^\infty\frac{1}{n!}X^n$.
Since $\eta\left(\sum_{n=0}^\infty\frac{1}{n!}X^n\right)=\zeta$, the condition (5) is a characterisation of Carlitz type of
exponential series.

The additive version of Theorem 1.3. is the following.

{\bf Theorem 1.4.} {\it Let $\sum_{n=0}^\infty a_n X^n\in\C[[X]]$ be a formal
power series with $a_1\not=0$. The following  statements are equivalent:

{\rm(1)} $a_0=0$ and $a_n=\frac{a_1}{(n-1)!}$, $\forall n\in\Z^+$;

{\rm(2)} $\eta\left(\sum_{n=0}^\infty a_n X^n\right)$ is an
additive arithmetical function;

{\rm (3)} $\sum_{n=0}^\infty a_nX^n\odot\left(\sum_{n=0}^\infty b_n X^n\cdot
\sum_{n=0}^\infty
c_nX^n\right)=$

\hspace{0.6cm} $=\left(\sum_{n=0}^\infty a_n X^n\odot\sum_{n=0}^\infty b_nX^n\right)\cdot
\sum_{n=0}^\infty c_n X^n+$

\hspace{0.6cm} $+\left(\sum_{n=0}^\infty a_nX^n\odot
\sum_{n=0}^\infty c_nX^n\right)\cdot\sum_{n=0}^\infty b_n X^n$,

\hspace{0.6cm}
$\forall
\sum_{n=0}^\infty b_n X^n$,
$\sum_{n=0}^\infty c_n X^n\in\C[[X]]$;

{\rm(4)} $\sum_{n=0}^\infty a_nX^n\odot\left(\sum_{n=0}^\infty
b_n X^n\cdot \sum_{n=0}^\infty b_n X^n\right)=$

\hspace{0.6cm} $=2\left(\sum_{n=0}^\infty a_nX^n\odot\sum_{n=0}^\infty
b_n X^n\right)\cdot\sum_{n=0}^\infty b_n X^n$,
$\forall\sum_{n=0}^\infty b_n X^n
\in\C[[X]]$;

{\rm(5)} $\sum_{n=0}^\infty 2^n a_n X^n=\sum_{n=0}^\infty \frac{2}{n!}X^n
\cdot\sum_{n=0}^\infty a_n X^n$.}

\begin{center}
{\bf 2.}
\end{center}

M\"obius categories were introduced by P. Leroux in [5]. In [6], P. Leroux has introduced the concept of M\"obius category of full binomial type and has denoted by $B(m)$
the parameters of a M\"obius category of full binomial type. $B(m)$ represent the total number
of decompositions into indecomposable factors of length 1 of a morphism of length $m$
$(\al=\al_1\al_2\dots\al_m$ with $\al_i\not=1$ is called a decomposition of degree $m$ of a morphism $\al$ and the
supremum of decomposition degree of $\al$ is called the length of $\al$ and is denoted by $\ell(\al))$.
In any M\"obius category of full binomial type, for two morphism $\al$ and $\beta$ with $\ell(\al)=\ell(\beta)=m$,
the following holds:
$$\left(\ba{c}
\al\\
k\ea\right)=\left(\ba{c}
\beta\\
k\ea\right)=\frac{B(m)}{B(k)\cdot B(m-k)}\buildrel{\mbox{not}}\over =
\left(\ba{c}m\\
k\ea\right)_\ell,$$
where $\left(\ba{c}
\al\\
k\ea\right)$ represent the number $|\{(\al',\al''):\al'\al''=\al, \ell(\al')=k\}|.$

Let $\cm$ be a M\"obius category of full binomial type and
$\ca_\cm$ the $\C$-algebra of arithmetical functions (with the
domain extended to $\Z^+\cup\{0\})$ under the usual pointwise
addition and scalar multiplication together with the
$\cm$-convolution defined by:
$$(f*_\cm g)(m)=\sum_{k=0}^m\left(\ba{c}
m\\
k\ea\right)_\ell f(k)g(m-k)\quad (m\in\Z^+\cup\{0\}).$$ Then
$\eta_\cm:\C[[X]]\to\ca_\cm$ defined by:
$$\eta_\cm\left(\sum_{n=0}^\infty a_n X^n\right)(m)=a_m B(m)\quad (m\in\Z^+\cup\{0\})$$
is an algebra isomorphism. We shall say that an arithmetical
function $f\in\ca_\cm$ is binomial multiplicative (binomial
additive) if $f(m+n)=f(m)\cdot f(n)$ $(f(m+n)=f(m)+f(n))$ for all
$m,n\in\Z^+\cup\{0\}$ and we have:

{\bf Theorem 2.1.} ([8]) {\it Let $\cm$ be a M\"obius category of full binomial type and
$\sum_{n=0}^\infty a_n X^n\in\C[[X]]$ with $a_1\not=0$. The following  statements are equivalent:

{\rm(1)} $a_n=\frac{a_1^n}{B(n)}$ for any $n\in\Z^+\cup\{0\}$;

{\rm(2)} $\eta_\cm\left(\sum_{n=0}^\infty a_n X^n\right)$ is a
binomial multiplicative arithmetical function;

{\rm (3)} $\sum_{n=0}^\infty a_nX^n\odot_\cm\left(\sum_{n=0}^\infty b_n X^n
\cdot\sum_{n=0}^\infty c_nX^n\right)=$

\hspace{0.6cm} $=\left(\sum_{n=0}^\infty a_n X^n\odot_\cm
\sum_{n=0}^\infty b_nX^n\right)\cdot \left(\sum_{n=0}^\infty a_n
X^n\odot_\cm\sum_{n=0}^\infty c_nX^n\right)$,

\hspace{0.6cm} $\forall \sum_{n=0}^\infty b_n X^n$,
$\sum_{n=0}^\infty c_n X^n\in\C[[X]]$;

{\rm(4)} $\sum_{n=0}^\infty a_nX^n\odot_\cm\left(\sum_{n=0}^\infty b_n X^n\cdot
\sum_{n=0}^\infty b_n X^n\right)=$

\hspace{0.6cm} $=\left(\sum_{n=0}^\infty a_nX^n\odot_\cm\sum_{n=0}^\infty
b_n X^n\right)\cdot\left(\sum_{n=0}^\infty a_n X^n\odot_\cm \sum_{n=0}^\infty b_n X^n\right)$,

\hspace{0.6cm} $\forall\sum_{n=0}^\infty b_n X^n
\in\C[[X]]$;

{\rm(5)} $\sum_{n=0}^\infty t(n) a_n X^n=\sum_{n=0}^\infty a_nX^n\cdot
\sum_{n=0}^\infty a_n X^n$,

\noindent where $t(n)=\sum_{k=0}^n\left(\ba{c}
n\\
k\ea\right)_\ell$ and the ``multiplication" $\odot_\cm$ of two formal power series
$\sum_{n=0}^\infty a_n X^n$ and $\sum_{n=0}^\infty b_n X^n$
is defined by:}
$$\sum_{n=0}^\infty a_n X^n\odot_\cm \sum_{n=0}^\infty b_n X^n=
\sum_{n=0}^\infty B(n) a_n b_n X^n.$$

Here (3) is also a characterization of Lambek type (distributivity over the product of series) and (5) is a
distributivity condition over a particular
product of series:
$\sum_{n=0}^\infty
\frac{1}{B(n)}X^n\cdot\sum_{n=0}^\infty\frac{1}{B(n)}X^n$.
Since $\eta_\cm\left(\sum_{n=0}^\infty\frac{1}{B(n)}X^n\right)=\zeta$,
the condition (5) is a characterisation of Carlitz type.

The additive version of Theorem 2.1. is the following.

{\bf Theorem 2.2.} {\it Let $\cm$ be a M\"obius category of full binomial type and
$\sum_{n=0}^\infty a_n X^n\in\C[[X]]$ with $a_1\not=0$.
The following  statements are equivalent:

{\rm(1)} $a_n=\frac{na_1}{B(n)}$ for any $n\in\Z^+\cup\{0\}$;

{\rm(2)} $\eta_\cm\left(\sum_{n=0}^\infty a_n X^n\right)$ is a
binomial additive arithmetical function;

{\rm (3)} $\sum_{n=0}^\infty a_nX^n\odot_\cm\left(\sum_{n=0}^\infty b_n X^n
\cdot\sum_{n=0}^\infty
c_nX^n\right)=$

\hspace{0.6cm} $=\left(\sum_{n=0}^\infty a_n X^n\odot_\cm\sum_{n=0}^\infty b_nX^n\right)\cdot
\sum_{n=0}^\infty c_n X^n+$

\hspace{0.6cm} $+\left(\sum_{n=0}^\infty a_nX^n\odot_\cm
\sum_{n=0}^\infty c_nX^n\right)\cdot\sum_{n=0}^\infty b_n X^n$,

\hspace{0.6cm} $\forall
\sum_{n=0}^\infty b_n X^n$,
$\sum_{n=0}^\infty c_n X^n\in\C[[X]]$;

{\rm(4)} $\sum_{n=0}^\infty a_nX^n\odot_\cm\left(\sum_{n=0}^\infty b_n X^n\cdot
\sum_{n=0}^\infty b_n X^n\right)=$

\hspace{0.6cm} $=2\cdot\left(\sum_{n=0}^\infty a_nX^n\odot_\cm\sum_{n=0}^\infty
b_n X^n\right)\cdot\sum_{n=0}^\infty b_n X^n$,

\hspace{0.6cm}
$\forall\sum_{n=0}^\infty b_n X^n
\in\C[[X]]$;

{\rm(5)} $\sum_{n=0}^\infty t(n) a_n X^n=2\cdot \sum_{n=0}^\infty \frac{1}{B(n)}X^n
\cdot\sum_{n=0}^\infty a_n X^n$.}

Straightforward verification shows that the characterization of the Lam\-bert-Carlitz
type of exponential series (Theorem 1.3.) is a special case of Theorem 2.1. if we find
a M\"obius category of full binomial type with the parameters
$B(n)=n!$. The category $\De$face
of injective increasing mappings of finite sets of the form $\{1,2,\dots, n\}$ with the initial
object $0$ is a M\"obius ca\-te\-gory of full
binomial type with the parameters $B(n)=n!$ (see [6]).
Then $\left(\ba{c}
n\\
k\ea\right)_\ell=\left(\ba{c}
n\\
k\ea\right)$, $t(n)=2^n$ and $\odot_{\De{\mbox{face}}}$ coincide with $\odot$.

Note that it is an embedding of $\De$face in the M\"obius
category of the reduced matrices over a finite field $GF(q)$
denoted by Red (see [6]). Red is a M\"obius category of full
binomial type with the parameters
$B(n)=[n]_q!=(1+q)(1+q+q^2)\dots(1+q+\dots+q^{n-1})$ ([6]). Thus
we obtain again the characterization of exponential series
(Theorem 1.3.) as a special case of the following
characterization of Lambek-Carlitz type over the M\"obius
category Red:

{\bf Theorem 2.3.} ([9]). {\it Let $\sum_{n=0}^\infty a_n X^n\in\C[[X]]$
be a formal power series with $a_1\not=0$.
The following  statements are equivalent:

{\rm(1)} $a_n=\frac{a_1^n}{[n]_q!}$ for any $n\in\Z^+\cup\{0\}$;

{\rm(2)} $\eta_{\mbox{\rm{\small{Red}}}}\left(\sum_{n=0}^\infty a_n X^n\right)$ is a
binomial multiplicative arithmetical function;

{\rm (3)} $\sum_{n=0}^\infty a_nX^n\odot_{\mbox{\rm{\small{Red}}}}\left(\sum_{n=0}^\infty b_n X^n
\cdot\sum_{n=0}^\infty
c_nX^n\right)=$

\hspace{0.6cm} $=\left(\sum_{n=0}^\infty a_n X^n\odot_{\mbox{\rm{\small{Red}}}}\sum_{n=0}^\infty b_nX^n\right)\cdot
\left(\sum_{n=0}^\infty a_nX^n\odot_{\mbox{\rm{\small{Red}}}}
\sum_{n=0}^\infty c_nX^n\right)$,

\hspace{0.6cm} $\forall
\sum_{n=0}^\infty b_n X^n$,
$\sum_{n=0}^\infty c_n X^n\in\C[[X]]$;

{\rm(4)} $\sum_{n=0}^\infty a_nX^n\odot_{\mbox{\rm{\small{Red}}}}\left(\sum_{n=0}^\infty b_n X^n\cdot
\sum_{n=0}^\infty b_n X^n\right)=$

\hspace{0.6cm} $=\left(\sum_{n=0}^\infty a_nX^n\odot_{\mbox{\rm{\small{Red}}}}\sum_{n=0}^\infty
b_n X^n\right)\cdot\left(\sum_{n=0}^\infty a_n X^n\odot_{\mbox{\rm{\small{Red}}}}
\sum_{n=0}^\infty b_n X^n\right)$,

\hspace{0.6cm} $\forall\sum_{n=0}^\infty b_n X^n
\in\C[[X]]$;

{\rm(5)} $\sum_{n=0}^\infty G_n(q) a_n X^n=\sum_{n=0}^\infty a_n X^n\cdot
\sum_{n=0}^\infty a_n X^n$,

\noindent where $G_n(q)$ are the Galois numbers and
$\sum_{n=0}^\infty a_n X^n\odot_{\mbox{\rm{\small{Red}}}} \sum_{n=0}^\infty b_n X^n=\sum_{n=0}^\infty
[n]_q!a_nb_n X^n$.}

The additive version of Theorem 2.3. follows from Theorem 2.2.

A survey of various special M\"obius categories of full binomial
type was given by P. Leroux [6]. All of them yields
characterizations of Lambek-Carlitz type of special formal power
series. The ``elementary'' M\"obius category of full binomial
type whose objects are the elements of $\Z^+\cup\{0\}$ and with
one morphism $(m,n):m\to n$ for each ordered pair $(m,n)$ in
which $m\le n$, induced the characterizations of Lambek-Carlitz
type of the formal power series $\sum_{n=0}^\infty a^n X^n$
(Theorem 2.1.) and $\sum_{n=0}^\infty na X^n$ (Theorem 2.2.).

{\bf Final Remark.} The isomorphisms between the algebra of formal
power series $\C[[X]]$ and some reduced incidence algebras of
M\"obius categories create relationships between the algebraic
(the Lambek type) and the combinatorial (the Carlitz type)
characterizations of formal power series. The
``number-theoretical component'' of the characterizations follows
from the fact that here the elements of the reduced incidence
algebras are arithmetical functions.

\noindent Address: Emil Daniel Schwab\\
The University of Texas at El Paso\\
Department of Mathematical Sciences\\
500, W. University Ave.\\
El Paso, Texas 79968, USA\\
e-mail: schwab@math.utep.edu

\end{document}